\documentclass[12pt]{amsart}
\usepackage{mathrsfs}
\usepackage{latexsym,amssymb,amsmath,amsfonts,graphicx,epstopdf,subfigure,float}
\usepackage[colorlinks=true,linkcolor=blue,citecolor=red,urlcolor=blue,]{hyperref}
\usepackage{pifont}
\usepackage{cite}
\usepackage{multirow}
\usepackage{appendix}
 \usepackage{color}

\newtheorem{theorem}{Theorem}[section]
\newtheorem{thm}[theorem]{Theorem}
\newtheorem{prop}{Proposition}[section]
\newtheorem{cor}[prop]{Corollary}
\newtheorem{lemma}[prop]{Lemma}
\newtheorem{remark}[prop]{Remark}
\newtheorem{defi}[prop]{Definition}

\newtheorem{example}[prop]{Example}
\numberwithin{equation}{section}

\def\cal{\mathcal }
\def\R{\mathbb R}

\def\Z{\mathbb Z}

\def\mathscr{\mathcal }

\newcommand{\ba}{{\mathbf{a}}}

\newcommand{\bd}{{\mathbf{d}}}

\newcommand{\bi}{{\mathbf{i}}}
\newcommand{\bj}{{\mathbf{j}}}

\newcommand{\SC}{{\mathcal C}}
\newcommand{\SD}{{\mathcal D}}

\newcommand{\SF}{{\mathcal F}}

\newcommand{\SP}{{\mathcal P}}

\newcommand{\diam}{\text{diam}}

\begin{document}

\title[]{Box-counting measure of metric spaces}


\author{Liang-yi Huang} \address{College of Mathematics and Statistics, Chongqing  University, Chongqing, 401331, China}
\email{liangyihuang@cqu.edu.cn}

\author{Hui Rao$\dagger$} \address{Department of Mathematics and Statistics, Central China Normal University, Wuhan, 430079, China
} \email{hrao@mail.ccnu.edu.cn
 }

 \author{Zhiying Wen} \address{Department of Mathematical Sciences, Tsinghua University, Beijing, 100084, China}
\email{wenzy@tsinghua.edu.cn}

 \author{Yan-li Xu}
\address{Department of Mathematics and Statistics, Central China Normal University, Wuhan,430079, China.}
\email{xu\_yl@mails.ccnu.edu.cn}

\date{\today}
\thanks {The work is supported by NSFS No.  11971195. }
\thanks{$\dagger$ Corresponding author.}


\begin{abstract} In this paper, we introduce a new notion called the \emph{box-counting measure} of a metric space.
We show that for a doubling metric space, an Ahlfors regular measure is always a box-counting measure; consequently, if $E$ is a self-similar set satisfying the open set condition, then the Hausdorff measure restricted to $E$
is a box-counting measure.
We show two classes of self-affine sets,   the generalized Lalley-Gatzouras type self-affine sponges and Bara\'nski carpets,  always admit box-counting measures; this also provides
 a very simple method to calculate the box-dimension of these fractals. Moreover, among others,  we show that if two doubling metric spaces admit box-counting measures, then  the multi-fractal spectra of the  box-counting measures coincide provided the two spaces are Lipschitz equivalent.
\end{abstract}
\maketitle


 \section{\textbf{Introduction}}

Let $X$ and $Y$ be two metric spaces with common Hausdorff dimension $s$.
Let $\mu$ and $\nu$ be the restrictions of the $s$-dimensional Hasudorff measure of $X$ and $Y$, respectively. It is well known that
if $\mu$ and $\nu$ are finite measures and $f:X\to Y$ is a bi-Lipschitz map, then
\begin{equation}\label{eq:key}
f^*\mu\sim \nu,
\end{equation}
where $f^*\mu(A):=\mu(f^{-1}(A))$; see for instance, Falconer \cite{Fal90}.
(Recall that two  measures $\mu$ and $\nu$ on $X$ are said to be \emph{equivalent}, and denoted by $\mu\sim \nu$,
if there exists $\zeta>0$ such that
$\zeta^{-1} \mu(\cdot) \leq \nu(\cdot)\leq \zeta \mu(\cdot).$)

Recall that  a measure $\mu$ on a metric space $X$ is said to be \emph{Ahlfors regular with index $s$} if
there is a constant $C>0$ such that $C^{-1}r^s\leq \mu(B(x,r))\leq Cr^s$
holds for any ball $B(x,r)$ with center $x\in X$ and radius $r\in (0,1)$.
 If  $\mu$ and $\nu$
are Ahlfors regular measures of $X$ and $Y$ with index $s$, respectively, then
 relation \eqref{eq:key} still holds. Indeed, in this case
$\mu\sim {\mathcal H}^s|_X$ and $\nu\sim {\mathcal H}^s|_Y$.
We ask the question that do there exist other metric spaces and measures such that \eqref{eq:key} holds?

Recently,  Rao \textit{et al.}  \cite{Rao2019}
found that \eqref{eq:key} holds for uniform Bernoulli measures on  Bedford-McMullen carpets.

\begin{prop}[\cite{Rao2019}] Let $X$ and $Y$ be two totally disconnected Bedford-McMullen carpets, let $\mu$ and $\nu$ be the uniform Bernoulli measures of $X$ and $Y$ respectively. If $f:X\to Y$ is bi-Lipschitz,
then $f^*\mu\sim \nu$.
\end{prop}

Consequently,  $\mu$ and $\nu$ have the same multi-fractal spectrum, and $\mu$ is doubling if and only if $\nu$ is doubling. Hence the multi-fractal spectrum and the doubling property can serve as new Lipschitz invariants. Thanks to these new Lipschitz invariants, Yang and Zhang \cite{Paper-3}
gave a complete Lipschitz
classification of totally disconnected Bedford-McMullen carpets whose Hausdorff dimension and box dimension coincide.

Recently,  Falconer,  Fraser and Kempton  \cite{FFK20}
introduced a new notion called  intermediate dimension.
Banaji and Kolossv\'ary \cite{B}
calculated the intermediate dimensions of Bedford-McMullen carpets, and   proved that
two  Bedford-McMullen  carpets have the same intermediate dimensions if and only if they have the same
multi-fractal spectrum. Therefore, they obtain that the following result without
assuming the totally disconnected condition.

 \begin{prop}[\cite{B}] Let $X$ and $Y$ be two    Bedford-McMullen carpets, let $\mu$ and $\nu$ be the uniform Bernoulli measures of $X$ and $Y$ respectively. If $f:X\to Y$ is bi-Lipschitz,
then $\mu$ and $\nu$ have the same multi-fractal spectrum.
\end{prop}

In the present paper, we introduce a new notion called the \emph{box-counting measure}, and we will show that \eqref{eq:key} holds for very general settings.

Let $(X, d_X)$ be a compact metric space and let $A\subset X$.
A family of balls with radius $\delta$ is called a \emph{$\delta$-ball-packing} of $A$, if they are disjoint and their centers are located
in $A$. For any $\delta>0$,   we define
\begin{equation*}
 \mathcal{N}_\delta(A):=\max \{\#\SP;~\text{$\SP$ is a $\delta$-ball-packing of  $A$}\},
\end{equation*}
where $\#{\mathcal P}$ denotes the cardinality of $\mathcal P$.
The upper and lower box dimensions of $X$ are defined by
$\overline{\dim}_B X=\overline{\lim}_{\delta\to 0} \frac{\log \mathcal{N}_\delta(X)}{-\log \delta}$
and $\underline{\dim}_B X=\underline{\lim}_{\delta\to 0} \frac{\log \mathcal{N}_\delta(X)}{-\log \delta}. $
If the two values coincide, then   the common value is called  the \emph{box dimension} of $X$
 and denoted   by $\dim_B X$.

 First, we define a new type  covering  of a metric space.

 \begin{defi}[Compact Vitalli-type covering] \emph{ Let $X$ be a compact metric space.
   For each $k\geq 1$, let
${\mathcal F}_k$ be a  covering of $X$. We call ${\mathcal F}=\{\SF_k\}_{k\geq 1}$
a \emph{compact Vitalli-type covering}  of $X$, if}

\emph{$(i)$ every element in $\SF_k$ is compact;}

\emph{$(ii)$  there is an integer $\lambda$ such that for all $k\geq 1$,
every point   $x\in X$ is covered by at most $\lambda$ elements of  $\SF_k$ ;}

\emph{$(iii)$  $\max \{\text{diam} A;~A\in \SF_k\}\to 0$ as $k\to \infty$.}
\end{defi}

Now we define the box-counting measure as follows.

\begin{defi}[Box-counting measure]\label{bcm}
\emph{ Let $X$ be a compact metric space such that  $\beta=\dim_B X$ exists.
Let $\mu$ be a finite Borel measure  on $X$.  Let ${\mathcal F}=\{\SF_k\}_{k\geq 1}$ be a
compact Vitalli-type  covering   of $X$.
We call $\mu$   an \emph{$\SF$-box-counting measure},
if   there is a constant $M>0$ such that
for any  $R\in \bigcup_{k\geq 1} \SF_k $,
\begin{equation}\label{eq:stable}
M^{-1}\mu(R)\leq \mathcal{N}_\delta(R)\delta^{\beta} \leq M \mu(R)
\end{equation}
holds for $\delta$ small enough (i.e., for  $0<\delta\leq\delta_0(R)$); in this case, we call
the triple $(X,\mu, \SF)$ a \emph{box-counting space}.
}
\end{defi}

\begin{example}\emph{Let $\epsilon>0$. Two points $x,y\in X$ are said to be \emph{$\epsilon$-equivalent} if there exists a sequence $\{x_1=x,x_2,\dots,x_{k-1},x_k=y\}\subset X$ such that $d_X(x_i, x_{i+1})\le\epsilon$ for $1\le i\le k-1$. An $\epsilon$-equivalent class is called an \emph{$\epsilon$-connected component of $X$}. Let $\SC(X, \epsilon)$ be the collection of all $\epsilon$-connected components of $X$.}

\emph{  If $X$ is a totally disconnected compact metric space, then
$$
\{\SF_k=\SC(X,  {1}/{2^k})\}_{k\geq 1}
$$
 is a compact Vitalli-type covering.}
\end{example}

For a large class of IFS's, their attractors admit a natural compact Vitalli-type covering.
An \emph{iterated function system} (IFS) is a  family of contractions $\Phi=\{\varphi_j\}_{j=1}^N$  on a compact set $X\subset \R^d$. In this paper, we will always assume that
   all $\varphi_j$'s are injections. The \emph{attractor} of the IFS is the unique nonempty compact set $E$ satisfying
$E=\bigcup_{j=1}^N\varphi_j(E)$; especially, it is called a \emph{self-similar set} if all $\varphi_j$'s are similitudes.
An IFS $\{\varphi_j\}_{j=1}^N$ is said to satisfy the \emph{open set condition} (OSC), if there is a bounded nonempty open set $U\subset\mathbb{R}^{d}$ such that for all $1\le i\le N$, $\varphi_i(U)\subset U$ and $\varphi_i(U)\cap\varphi_j(U)=\emptyset$ for $1\le i\ne j\le N$.
See \cite{Hut81, Fal90}.

Denote $\Sigma=\{1,\dots, N\}$ and $\Sigma^*=\bigcup_{n\geq 0} \Sigma^n$.
For $I=i_1\dots i_n\in \Sigma^*$, we call
$E_I=\varphi_I(E)=\varphi_{i_1}\circ\dots\circ\varphi_{i_n}(E)$
an $n$-th \emph{cylinder} of $E$.
Let $x\in E$, we call $J=(j_k)_{k\geq 1}\in \Sigma^\infty$ a \emph{coding}
of $x$ if $x\in \varphi_{j_1\dots j_n}(E)$ for all $n\geq 1$.

\begin{defi}\emph{
Let $\SF_k$ be the collection of $k$-th cylinders of an IFS $\Phi$.
Clearly $\{\SF_k\}_{k\geq 1}$ is a compact Vitalli-type covering of the attractor $E$
if and only if  the number of codings of  points in $E$ is uniformly bounded.
  Let $\mu$ be a finite  Borel measure  on $E$.
We call $\mu$ a \emph{cylinder box-counting measure} if  $(E, \mu, \SF=\{\SF_k\}_{k\geq 1})$ is a
box-counting space.
}\end{defi}

\begin{figure}[H]
  \centering
  \includegraphics[width=4 cm]{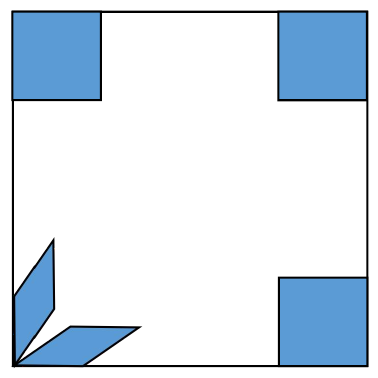}
\caption{Let $F$ be the attractor of the above self-affine IFS, then $F$ satisfies the
OSC, but the left-bottom point has infinitely many codings.}
\label{Q'}
\end{figure}

A metric space $X$ is said to be a \emph{doubling space}, if there is a constant $C>0$ such that for any $x\in X$ and $r>0$, $B(x,2r)$ can be covered by    $C$ numbers of  balls of radius $r$ (see \cite{Heinonen}). For a doubling space, the box-counting measure is a generalization of Ahlfors regular measure.

\begin{thm}\label{thm:regular}
Let $X$ be a compact doubling   space. Then every Ahlfors regular measure of $X$  is a   box-counting measure.
\end{thm}

\begin{example} \label{thm:SSS}\emph{
Let $E\subset \R^d$ be a self-similar set satisfying the OSC. Let $s=\dim_H E$.
Denote ${\mathcal H}^s$ the $s$-dimensional Hausdorff measure.
It is well-known that ${\mathcal H}^s|_E$  is an Ahlfors measure, and hence it is also a
cylinder box-counting measure. A direct proof is given in Theorem \ref{thm:OSC}.}
\end{example}

It is well known that on a metric space, any two Ahlfors regular measures are equivalent. Similarly, we have

\begin{thm}\label{thm:main-1}
\emph{Let $(X, \mu, \SF)$ and  $(X, \mu', \SF')$ be two box-counting spaces.  Then $\mu\sim\mu'$.}
\end{thm}


 Next, we investigate Lipschitz invariants related to box-counting measures.

\begin{thm}\label{thm:main-2}  Let $X$ and $Y$ be two compact doubling spaces.
 Suppose $(X, \mu, \SF)$ is a  box-counting space, and  $f:~X\to Y$ is bi-Lipschitz. Then
 $(Y, f^*\mu, f(\SF))$ is a   box-counting  space.
\end{thm}

As a corollary of  Theorem \ref{thm:main-1} and \ref{thm:main-2}, we have

\begin{cor}\label{cor:main}  Let $X$ and $Y$ be two  compact doubling spaces.
 Suppose $\mu$ and $\nu$ are  box-counting measures
of $X$ and $Y$,  respectively.
 If $f:~X\to Y$ is bi-Lipschitz, then
 $$f^*\mu\sim \nu.$$
Consequently,  $\mu$ and $\nu$ have the same multi-fractal spectrum, and $\mu$ is doubling if and only if $\nu$ is doubling.
\end{cor}

\begin{proof} Since both $(Y, f^*\mu, f(\SF))$ and $(Y, \nu, \SF')$ are box-counting spaces, by
Theorem \ref{thm:main-1}, we have $f^*\mu\sim \nu$. Therefore, $\mu, f^*\mu$ and $\nu$ have
the same multi-fractal spectrum and the same doubling property.
\end{proof}

We remark that the above corollary generalizes the corresponding results in
Rao, Yang and Zhang \cite{Rao2019} and Banaji and Kolossv\'ary \cite{B}.

Finally, we show that several class of diagonal self-affine sets admits  cylinder box-counting measures.

\begin{defi}\emph{ We call $f:\R^d\to \R^d$, $f(x)=Tx+b$ a \emph{diagonal self-affine mapping}
if $T$ is a $d\times d$ diagonal matrix such that all the diagonal entries are positive numbers.
An IFS $\Phi=\{\phi_j(x)\}_{j=1}^m$ is called a \emph{diagonal self-affine IFS}
if all the maps $\phi_j(x)$ are  diagonal self-affine contractions; the attractor is called
a \emph{diagonal self-affine sponge},  which we denote  by $\Lambda_\Phi$. Without loss of generality, we will always assume that  $\Lambda_\Phi\subset [0,1]^d$. }
\end{defi}

Das and Simmons  \cite{Das16} gave an equivalent definition as following.
 For each $i\in \{1,\ldots,d\}$, let $A_i$ be a finite index set, and let $\Phi_i=\{\phi_{a,i}\}_{a\in A_i}$ be a collection of contracting similarities of $[0,1]$ coming from the $i$-th components
 of maps in $\Phi$, called
the \emph{base IFS in coordinate} $i$. Let $A =A_1\times\cdots \times A_d$, and for each $\ba=(a_1,\ldots,a_d)\in A,$ define the contracting affine map
$\phi_\ba: [0,1]^d \rightarrow [0,1]^d$ by
$$
\phi_\ba(x_1,\ldots,x_d)=(\phi_{\ba,1}(x_1),\dots,\phi_{\ba,d}(x_d)),
$$
where $\phi_{\textbf{a},i}$ is shorthand for $\phi_{a_i,i}$ in the formula above.
Given $\SD\subset A$, then the collection
 $\Phi=\{\phi_\ba\}_{\ba\in \SD}$ is a diagonal self-affine IFS.

%

\begin{remark}
\emph{Recently, there are a lot of works on diagonal self-affine sponges,   Feng and Wang \cite{FengWang05}, Bara\'nski \cite{Baranski07},  Mackay \cite{MM11}, Fraser \cite{Fra13},
Das and Simmons \cite{Das16},   Banaji and Kolossv\'ary \cite{B}  on dimensions;
King  \cite{King95}, Jordan and Rams \cite{JR11}, Olsen \cite{Olsen07}, Reeve \cite{Reeve10} on multi-fractal formalism; Li, Li and Miao \cite{Miao2013}, Miao, Xi and Xiong \cite{Miao2017},
Liang, Miao and Ruan \cite{LMR22},
Rao, Yang and Zhang \cite{Rao2019} on metric and topology classifications.}
\end{remark}


\begin{figure}[H]
  \centering
  \includegraphics[width=12 cm]{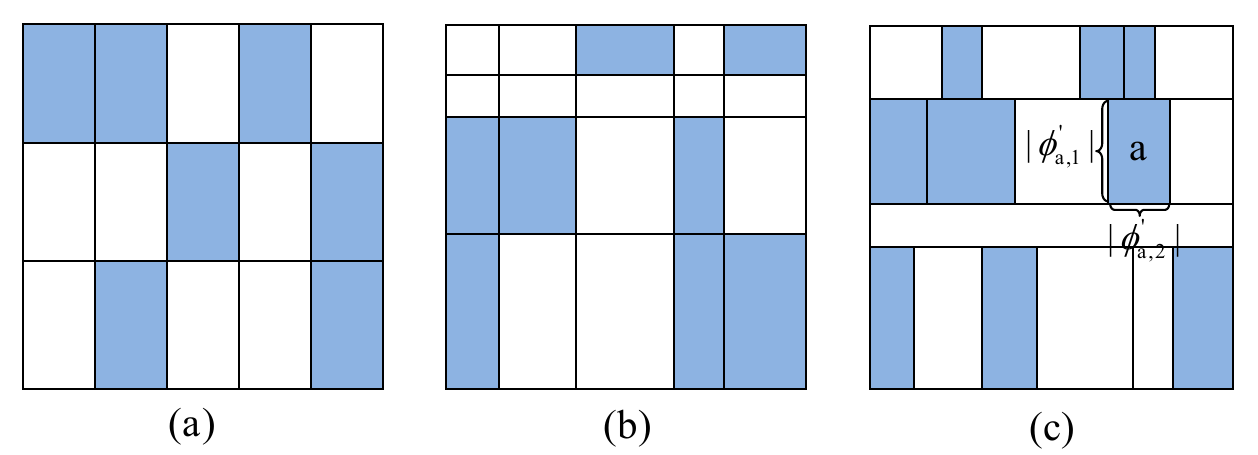}\\
  \caption{(a) Bedford-McMullen carpet; (b) Bara$\acute{\text{n}}$ski carpet; (c) Lalley-Gatzouras carpet.}\label{fig1}
\end{figure}

Now let $\Lambda_{\Phi}$ be a diagonal self-affine sponge.
Given a permutation $\tau$ of the coordinate set $\{1,\dots, d\}$,
a sequence $\{\beta_j\}_{j=1}^d\subset [0,1]$    can be defined inductively by \eqref{eq:flag}; see Section 4 for details.
For $\ba\in \SD$, we define
\begin{equation}\label{mu-p}
p_\ba=\prod_{j=1}^d (\phi_{\ba,_j}')^{\beta_j},
\end{equation}
 which is a probability weight.
Denote $\beta(\tau)=\sum_{j=1}^d \beta_j$ and let $\mu_\tau$ be the Bernoulli measure determined
 by the above probability weight.  This construction has been given in \cite{Baranski07} and \cite{Lalley92}
 for calculating the box dimensions of self-affine sets.

 A diagonal self-affine sponge $\Lambda_{\Phi}$ is said to satisfy the \emph{coordinate ordering condition}
 if there is a permutation $\tau$  such that $\{\phi_{\ba,j}'\}_{j=1}^d$
is a strictly decreasing sequence for each $\ba\in \SD$; is said to satisfy the   \emph{weak coordinate ordering condition} if $\{\phi_{\ba,j}'\}_{j=1}^d$ is non-increasing.

  A diagonal self-affine sponge is said to be of the \emph{generalized Lalley-Gatzouras type} (resp. \emph{Lalley-Gatzouras type}), if
it satisfies the weak coordinate ordering condition (resp. coordinate ordering condition)
as well as the neat projection condition.
 (See Section 4 for details).

\begin{thm}\label{thm:sponge} Let $\Lambda_{\Phi}$ be a diagonal self-affine sponge  of generalized Lalley-Gatzouras type. Let $\tau$ be the permutation such that $\{\phi_{\ba,j}'\}_{j=1}^d$
is a non-increasing sequence for each $\ba\in \SD$.   Then
$\dim_B \Lambda_{\Phi}=\beta(\tau)$,
 and the Bernoulli measure $\mu_\tau$ defined
  by the weight in \eqref{mu-p} is a cylinder box-counting measure.
\end{thm}

\begin{remark}\emph{Actually, we will show that for a cylinder $R=\Lambda_I$,
$\delta_0(R)$ can be chosen to be the length of the shortest side of $\phi_I([0,1]^d)$.}
\end{remark}

Let $\Lambda$ be a Bara\'nski carpet (see Section 5 for details).
Let $\tau_1$ be the identity map  and $\tau_2$ be the permutation $1\mapsto 2, 2\mapsto 1$.
 Following \cite{Baranski07},  denote  $$\alpha=\beta_{\tau_1}, \mu_A=\mu_{\tau_1} \text{ and }
   \beta=\beta_{\tau_2} ,  \mu_B=\mu_{\tau_2}.$$
Bara\'nski \cite{Baranski07} proved that
  $\dim_B \Lambda=\max\{\alpha,\beta\}.$

\begin{thm}\label{thm:Baranski}
Let $\Lambda$ be a  Bara\'nski carpet.
 Then
$$
\mu=\left \{
\begin{array}{ll}
\mu_A, & \text{ if }\alpha>\beta;\\
\mu_B, & \text{ if }\alpha<\beta;\\
\mu_A+\mu_B, &\text{ if }\alpha=\beta
\end{array}
\right .
$$
 is   a cylinder box-counting measure.
\end{thm}

\begin{remark}
\emph{ In the proof of  Theorem  \ref{thm:sponge} and \ref{thm:Baranski}, as a by-product, we  obtain
the box dimensions of the generalized Lalley-Gatzouras type sponges and the Bara\'nski carpets, respectively.
The first one is also calculated in Lalley and Gatzouras \cite{Lalley92} and  Kolossv\'ary \cite{K21}.
Due to the notion of box-counting measure, our proofs are considerably simpler than
 \cite{Lalley92}, \cite{K21} and  Bara\'nski
\cite{Baranski07}.
}
\end{remark}

\begin{remark} \emph{The box dimension is   closely related to
  Minkowski content (see Section \ref{sec:Minkowski} for precise definition).
 For applications of Minkowski content in analysis,  we refer to  Lapidus and   Pomerance \cite{LP93}, Falconer \cite{Fal95} and the references therein.}
\emph{ We show that if $(X,\mu, \SF)$ is a box-counting space, denote $\beta=\dim_B X$,
  then both the lower and upper $\beta$-dimensional Minkowski contents are finite
  (Lemma \ref{lem:content}).  }
  \end{remark}

 \begin{example}[Self-similar sets admitting no box-counting measures]\label{exam:four}
\emph{
Let $F_\lambda$ be the  self-similar set generated by the IFS
$$\left \{f_1(x)=\frac{x}{3}, \ f_2(x)=\frac{x+1}{3}, \ f_3(x)=\frac{x+\lambda}{3} \right \}.$$
When $\lambda$ is an irrational number, it was shown by Kenyon \cite{Kenyon97} that the Lebesgue measure of $F_\lambda$
 is zero, and it was proved by Hochman \cite{Hochman12} that $\dim_H F_\lambda=\dim_B F_\lambda=1$.
 Therefore the Minkowski content of $F_\lambda$ is zero, and $F_\lambda$ does not admit any box-counting measure
 by Lemma \ref{lem:content}.
}
\end{example}

 \begin{remark}\emph{
 Another motivation of this paper is to investigate the  \emph{locally measure preserving property} of bi-Lipschitz maps.
  Cooper and Pignataro \cite{CP88}, Xi and Ruan \cite {XR2008} proved that if $E$ and $F$ are two self-similar sets
 satisfying the strong separation condition,
 and $f:E\to F$ is bi-Lipschitz, then there is a cylinder $U$ of $E$ such that
 $f: U\to f(U)$ preserves the Hausdorff measure in dimension $\dim_H E$.
 This   property   plays an important r\^ole
 in the Lipschitz classification of self-similar sets, see \cite{FaMa92, RRW12, RZ15}, etc.}

\emph{In a sequential paper,    Yang and Zhang \cite{Yang22}  proved that if $(E,\mu)$ and $(F,\nu)$ are two box-counting spaces,
then the above measure preserving property is valid provided that}
 \emph{$(i)$ both $E$ and $F$ are \emph{perfectly disconnected};  $(ii)$ both $\mu$ and $\nu$ are \emph{arithmetically doubling}.
For precise definition of these teminologies, we refer to \cite{Yang22}.}
 \end{remark}

 Motivated by Banaji and Kolossv\'ary \cite{B}, we ask the following question.

 \medskip

\noindent \textbf{Open probelm 1.} Let $\Lambda_1$ and $\Lambda_2$ be two self-affine sponges
of generalized Lalley-Gatzouras type (or Bara\'nski carpets)  and let $\mu_1$ and $\mu_2$ be their
box-counting measures, respectively. Is it true that $\Lambda_i$ $(i=1,2)$ have the same
intermediate dimensions if and only if $\mu_i$ $(i=1,2)$ have the same multi-fractal
spectrum? (It is well-known that the calculation of intermediate dimensions and multi-fractal spectra are both tedious for self-affine sets.)

\medskip

The paper is organized as follows.  In Section \ref{sec:proof-thm-1}, we prove Theorem \ref{thm:regular}.  Theorem \ref{thm:main-1} and Theorem \ref{thm:main-2}
are proved in Section \ref{sec:property}.
In Section \ref{sec:sponge}, we study the box-counting measure of   self-affine sponges of generalized Lalley-Gatzouras type. In Section \ref{sec:Baranski}, we prove Theorem \ref{thm:Baranski}. In Section \ref{sec:symbolic}, we discuss the box-counting measures of some symbolic spaces.
The relation between box-counting measure and Minkowski content is discussed in Section \ref{sec:Minkowski}.

 \section{\textbf{Ahlfors regular measures are box-counting measures}}\label{sec:proof-thm-1}

Let $X$ be a metric space admitting an Ahlfors regular measure $\mu$. It is  well-known  that
   $\dim_H X=\dim_B X=s$.
Clearly, if $\mu'$ is another Ahlfors regular measure of $X$, then $\mu\sim \mu'$.

  For a set $A\subset X$, let
$$[A]_\delta:=\{x\in X;~d_X(x,a)\le\delta\text{ for some }a\in A\},$$
 and we call it the \emph{$\delta$-neighbourhood} of $A$. Let $\mathcal{L}$ be the collection of all compact subsets of $X$. If $A, B\in\mathcal{L}$, the \emph{Hausdorff metric} $d_H$ is defined as follows,
$$
d_H(A,B)=\inf\{\delta;~ A\subset [B]_\delta\text{ and }B\subset [A]_\delta\}.
$$
It is well-known that the space $\mathcal{L}$ is compact in the Hausdorff metric $d_H$ (see for instance,
\cite[Sec. 45]{Munkres}).

Through the whole paper, when we write $f(x)\asymp g(x),~x\in X$, we mean that there is a constant
$C$ independent of $x$  such that
$C^{-1}f(x)\leq g(x)\leq Cf(x)$ for all $x\in X$. In this case, we say that
$f(x)$ and $g(x)$ are \emph{comparable} for $x\in X$.

\begin{proof}[\textbf{Proof of Theorem \ref{thm:regular}}]
Let $X$ be a compact doubling space with $\beta=\dim_BX$, and let $\mu$ be an Ahlfors regular measure on $X$.

 First, we construct a compact Vitalli-type covering  of $X$.
That $X$ is compact implies $X$ is totally bounded.
Let $\{r_k\}_{k\geq 1}$ be a sequence of reals decreasing to $0$.  For each integer $k\geq 1$,
there exists a finite set $A'_k\subset X$ such that
$$\SF_k'=\{B(a,r_k);~a\in A'_k\}$$
is a covering of $X$. By the $5r$-covering lemma,  we can choose $A_k\subset A'_k$ such that
the balls $B(a, r_k), a\in A_k$, are disjoint and
$$\SF_k=\{B(a,5r_k);~a\in A_k\}$$
is a covering of $X$.

Pick $x\in X$. Suppose $x$ is covered by $p$ numbers of balls in $\SF_k$.
Then the ball $B(x, 10r_k)$ contains at least $p$ numbers of balls in $\SF_k'$.
Since $X$ is a doubling space, we conclude that $p\leq \lambda=C^4$ where $C$
is the constant in the definition of doubling space.   This proves that $\{\SF_k\}_{k\geq 1}$
is a compact Vitalli-type covering.

Next, let $\mathcal{N}'_\delta(R)$ be the minimal cardinality of $\delta$-ball-coverings of $R$.
 By the doubling property of $X$,  we have  that for any $R\subset X$,
 $$
 \mathcal{N}_\delta(R)\geq \mathcal{N}'_{2\delta}(R)\geq  C^{-1}\mathcal{N}'_{\delta}(R).
 $$

Finally, we show that for any $R\in \bigcup_{k=1}^\infty\SF_k$, $\mu(R)$ and $\mathcal{N}_\delta(R)\delta^{\beta}$ are comparable for $\delta$ small enough.
 Since $\mu$ is an Ahlfors regular measure,
 for any $x\in X$ and $\delta\in (0,1)$,
$$C_1^{-1} \delta^{\beta}\leq \mu(B(x,\delta))\leq C_1\delta^{\beta}$$
 for some $C_1>0$.
On one hand,  let $\{B_i\}_{i=1}^k$ be a minimal $\delta$-ball-covering of $R$, we deduce that
$$\mu(R)\leq \sum_{i=1}^k \mu(B_i)\leq \mathcal{N}'_\delta(R)C_1\delta^{\beta}\le CC_1\mathcal{N}_\delta(R)\delta^{\beta}.$$
 On the other hand, let $\{B'_j\}_{j=1}^\ell$ be a maximal $\delta$-ball-packing of $R$,  we have
$$
\mu([R]_\delta)\geq \sum_{j=1}^\ell \mu(B'_j)\geq \mathcal{N}_\delta(R) C_1^{-1}\delta^\beta.
$$
Note that $2\mu(R)\geq \mu([R]_\delta)$ holds for $\delta$ small, we obtain the other side estimation and finish the proof.
\end{proof}

For a self-similar set $E$ satisfying the OSC, the Hausdorff dimension $s$
is given by Moran's formula $\sum_{j=1}^N r_j^s=1$, where $r_j$'s are contraction ratios
of maps in the IFS. We call the Bernoulli measure
with probability weight $(r_j^s)_{j=1}^N$ the \emph{canonical Bernoulli measure} of $E$.

\begin{thm}\label{thm:OSC} Let $E\subset \R^d$ be a self-similar set satisfying the OSC.
Let $\mu$ be the canonical Bernoulli measure of $E$. Then there is a constant
$M>0$ such that for any cylinder $E_I$ and \emph{$\delta\le \diam(E_I)$},  it holds that
$$
M^{-1}\mu(E_I)\delta^{-\beta}\leq {\mathcal N}_\delta(E_I)\leq M\mu(E_I)\delta^{-\beta}.
$$
\end{thm}

\begin{proof} Let $\{\varphi_j\}_{j=1}^N$ be an IFS generating $E$.
Denote $\beta=\dim_B E$.
For $J\in \{1,\dots, N\}^*$,
let $c_J$ be the contraction ratio of $\varphi_J$. Let $\delta\in (0,\diam~ E]$. We set
$$
\Omega_\delta=\{J=j_1\dots j_k;~ c_J\leq \delta<c_{j_1\dots j_{k-1}}\}.
$$
It is well known that (see, for instance, Falconer \cite{Fal90})
$$
{\mathcal N}_\delta(E)\asymp \# \Omega_\delta \asymp \delta^{-\beta}.
$$
Therefore, since $\varphi_I$ is a similarity, we obtain that for $\delta\in(0, \diam~E_I]$,
$$
{\mathcal N}_\delta(E_I)={\mathcal N}_{\delta/c_I}(E)  \asymp
 (\delta/c_I)^{-\beta}=\mu(E_I)\delta^{-\beta}.
$$
The theorem is proved.
\end{proof}

 The above theorem shows that $\delta_0(R)$ in Definition \ref{bcm} can be chosen to be $\diam (R)$.

\section{\textbf{Proofs of Theorem \ref{thm:main-1} and Theorem \ref{thm:main-2}}}\label{sec:property}

In this section, we  investigate the Lipschitz invariants related to   box-counting spaces.

\subsection{\textbf{Equivalence of box-counting measures}}

\begin{proof}[\textbf{Proof of Theorem \ref{thm:main-1}}]
Let $\SF=\{\SF_k\}_{k\ge 1}$ and $\SF'=\{\SF'_k\}_{k\ge 1}$ be the compact Vitalli-type coverings in
the theorem.
Let $\lambda$ be the maximal covering
multiplicity of $\SF_k$ as well as $\SF_k'$ for all $k\geq 1$.
For $R\in \bigcup_{k=1}^\infty\SF_k$, denote
$$
R_k=\bigcup_{A'\in \SF'_k \text{ and }A'\cap R\neq \emptyset}A'.
$$
First, since the maximal diameter of elements of $\SF'_k$ tends to $0$ as $k\to \infty$,
we conclude that $d_H(R_k, R)\to 0$ as $k\to \infty$.
Secondly, since $R$ is compact, we have  that $R_k$   decreases to $R$ and it follows that
$\lim_{k\to \infty} \mu'(R_k)=\mu'(R).$

Since $(X,\mu,\SF)$ and $(X,\mu',\SF')$ are box-counting spaces, there are positive constants $M$ and $M'$ such that
\begin{align*}
M^{-1}\mu(A)&\leq \mathcal{N}_\delta(A)\delta^{\beta} \leq M \mu(A)\text{ for any } A\in\SF_k, k\ge 1;\\
M'^{-1}\mu(A')&\leq \mathcal{N}_\delta(A')\delta^{\beta} \leq M' \mu(A')\text{ for any } A'\in\SF'_k, k\ge 1.
\end{align*}
Choose $k$ large enough so that
$\mu'(R_k)<2\mu'(R)$, then we have
$$
\begin{array}{rl}
\mathcal{N}_\delta(R)\leq \mathcal{N}_\delta(R_k) &\leq  \sum\limits_{A'\subset R_k} \mathcal{N}_\delta(A')\\
&\leq \sum\limits_{A'\subset R_k} M'\mu'(A')\delta^{-\beta}\\
&\leq \lambda M'\mu'(R_k)\delta^{-\beta}\\
&\leq 2 \lambda M'\mu'(R)\delta^{-\beta}.
\end{array}
$$
From the above relations we conclude that $\mu(R)\leq 2\lambda MM'\mu'(R)$.

Since $\bigcup_{k=1}^\infty\SF_k$ generates the Borel algebra of $X$, we conclude that
$$\mu(A)\leq 2\lambda MM'\mu'(A)$$
for any Borel set $A\subset X$.
By symmetry, we also have $\mu'(A)\leq 2\lambda MM'\mu(A)$. Hence $\mu\sim \mu'$.
The theorem is proved.
\end{proof}

\subsection{\textbf{Bi-Lipschitz maps between doubling spaces}}\label{sec:proof-thm-main}
\

 Two metric spaces $(X, d_X)$ and $(Y, d_Y)$ are said to be \emph{Lipschitz equivalent},
and denoted by $(X,d_X)\sim (Y,d_Y)$, if there exist a bijection
$f:~X\rightarrow Y$ and a constant $C>0$ such that
$$C^{-1}d_X(x,y)\leq d_Y(f(x),f(y))\leq C d_X(x,y), \ \text{ for all } x,y\in X;$$
 in this case, we call $f$ a \emph{bi-Lipschitz map}.


\begin{lemma}\label{lem:box-counting}
Let $X$ and $Y$ be two  doubling  spaces, and let $f:~X\to Y$ be a bi-Lipschitz map. Then there exists a constant $C>0$ such that for any $\delta>0$,
\begin{equation}\label{eq-box-counting}
C^{-1}\mathcal{N}_\delta(X) \leq \mathcal{N}_\delta(Y)\leq C \mathcal{N}_\delta(X).
\end{equation}
\end{lemma}

\begin{proof}
Let $C_1>0$ be a  Lipschitz constant of $f$. For $r>0$ and  $x\in X$, we have
$$
B(f(x),C_1^{-1}r)\subset f(B(x,r)),
$$
it follows that $\mathcal{N}_\delta(X)\leq \mathcal{N}_{C_1^{-1}\delta}(Y)$.

Since $Y$ is a doubling space, there exists a constant $C_2>0$ such that
 for any $y\in Y$ and $r>0$, $B(y, r)$ contains at most $C_2$ disjoint balls of
 radius $C_1^{-1}r/3$.
Let $\SP_1=\{B(x_1, C_1^{-1}\delta), \dots, B(x_k, C_1^{-1}\delta)\}$
 be a $(C_1^{-1}\delta)$-ball-packing of $Y$ such that $k=\mathcal{N}_{C_1^{-1}\delta}(Y)$,
let $\SP_2=\{B(y_1, \delta), \dots, B(y_\ell, \delta)\}$
be a $\delta$-ball-packing of $Y$   such that $\ell=\mathcal{N}_{\delta}(Y)$.
 Then for each $x_i (1\le i\le k)$, there exists $y_j (1\le j\le \ell)$
 such that  $d_Y(x_i, y_j)\le 2\delta$.
 Let $\SP_3=\{B(y_1, 3\delta), \dots, B(y_\ell, 3\delta)\}$,
 we see that
  $\bigcup\SP_{1}\subset\bigcup\SP_{3}$. Therefore,
$$
\mathcal{N}_{C_1^{-1}\delta}(Y)=k \leq C_2\ell =C_2\mathcal{N}_\delta(Y).
$$
This proves the first inequality of \eqref{eq-box-counting}. By symmetry, we have the second inequality.
\end{proof}

\begin{proof}[\textbf{Proof of Theorem \ref{thm:main-2}.}]
 First, we claim  $f(\SF)$ is a compact Vitalli-type covering   of $Y$.
  Clearly, $f(\SF_k)$ is a covering of $Y$  preserving the maximal covering multiplicity $\lambda$, the maximal  diameter (up to the bi-Lipschitz constant $C$), and  the compactness.
Our claim is proved.

Let $\beta$ be the common value of the box dimension of $X$ and $Y$.
Let $R\in \bigcup_{k=1}^\infty \SF_k$.
 By Lemma  \ref{lem:box-counting} we have
\begin{equation}\label{eq:est_6}
\mathcal{N}_\delta(f(R))\asymp \mathcal{N}_\delta (R).
\end{equation}
  That $\mu$  is a box-counting measure of $X$ implies for $\delta\leq \delta_0(R)$,
\begin{equation}\label{est_4}
\mathcal{N}_\delta(R)\asymp  \mu(R)\delta^{-\beta}=(f^*\mu)(f(R))\delta^{-\beta}.
\end{equation}
Combining  \eqref{eq:est_6} and \eqref{est_4}, we obtain that
$$ \mathcal{N}_\delta(f(R))\asymp (f^*\mu)(f(R))\delta^{-\beta},\quad \text{ for }\delta\leq \delta_0(R).
$$
Hence $(Y, f^*\mu, f(\SF))$ is a box-counting measure  space.
 The theorem is proved.
\end{proof}

\section{\textbf{Box-counting measures of  self-affine sponges of generalized Lalley-Gatzouras type}}\label{sec:sponge}

Let $\Phi=\Phi_\SD$ be a diagonal self-affine IFS on $\mathbb{R}^d$ and let
$\Lambda_\Phi$ be the attractor of $\Phi$.
That $\Phi$ is said to satisfy the \emph{weak coordinate ordering condition} if there exists a permutation $\sigma$ of $\{1,2,\dots,d\}$ such that
$$
\phi'_{\ba,\sigma(1)} \geq \dots\geq \phi'_{\ba,\sigma(d)}, \quad \text{ for all } \ba\in\SD.
$$
Without loss of generality,   we always assume that
\begin{equation}\label{order}
\phi'_{\ba,1} \geq \dots\geq \phi'_{\ba,d}, \quad  \text{ for all }\ba\in\SD.
\end{equation}

Let $\pi_j: \mathbb{R}^d\rightarrow\mathbb{R}^j$ be the projection   $\pi_j(x_1,\dots,x_d)=(x_1,\dots,x_j).$
We define the $j$-th projection IFS of $\Phi$ to be
\begin{equation*}
\Phi_{\{1,\dots, j\}}=\left\{(\phi_{\textbf{d},1},\dots, \phi_{{\mathbf d},j})\right\}_{\textbf{d}\in{\pi_{{j}}}(\SD)},
\end{equation*}
which is a diagonal self-affine IFS on $\R^j$. (Here we emphasize that each map occurs at most once in the above IFS.)
Clearly,
\begin{equation*}
\Lambda_j:=\pi_j(\Lambda_\Phi)
\end{equation*}
is the attractor of the IFS $\Phi_{\{1,\dots, j\}}$.


\begin{defi} [\cite{Das16}] \label{def:good}\emph{
Let $\Lambda_\Phi$ be a diagonal self-affine sponge satisfying \eqref{order}.
 We say $\Phi$ satisfies the \emph{neat projection condition},
 if for each $j\in \{1,\dots, d\}$, the IFS $\Phi_{\{1,\dots,j\}}$ satisfies
 the OSC with the open set $\mathbb{I}^j=(0,1)^j$,  that is,
$$
\left\{\phi_{\bd,\{1,\dots, j\}}(\mathbb{I}^j)\right\}_{\bd\in\pi_j(\SD)}
$$
are disjoint.}
\end{defi}

In the rest of this section, we always assume that  $\Lambda_\Phi$   satisfies the
weak coordinate ordering condition \label{eq:order} as well as the neat projection condition, in other words,
$\Lambda_\Phi$ is of  \emph{generalized Lalley-Gatzouras type}.

Now we define a sequence of Bernoulli measures related to $\Lambda_j, j=1,\dots,d$.

First, we define a sequence $\{\beta_j\}_{j=1}^d$
related to $\Lambda_\Phi$.
Let $\beta_1>0$ be the unique real number satisfying
$$\sum\limits_{f_1\in \Phi_{\{1\}}} (f_1')^{\beta_1}=1.$$
If $\beta_1,\dots, \beta_{j-1}$ are defined, we define $\beta_j>0$ to be the unique real number such that
\begin{equation}\label{eq:flag}
\sum\limits_{(f_1,\dots, f_j)\in \Phi_{\{1,\dots, j\}}} \prod_{k=1}^j (f_k')^{\beta_k}=1.
\end{equation}

Next, for ${\mathbf f}=(f_1,\dots, f_j)\in \Phi_{\{1,\dots, j\}}$, define
$$
p_{\mathbf f}=\prod_{k=1}^j (f_k')^{\beta_k}.
$$
Let $\mu_j$ be the Bernoulli measure on $\Lambda_j$ defined by the weight
$(p_{\mathbf f})_{{\mathbf f} \in \Phi_{\{1,\dots, j\}}}$.
Especially, $\mu_d$ is the measure  in Theorem \ref{thm:sponge}.

\indent For $I=i_1\dots i_n\in\SD^n$, we call
$\phi_I(\Lambda_\Phi)$ an \emph{$n$-th cylinder} of $\Lambda_\Phi$, and call $\phi_{i_1\dots i_{n-1}}(\Lambda_\Phi)$ the \emph{ancestor} of $\phi_I(\Lambda_\Phi)$.
For a cylinder $R=\phi_I(\Lambda_\Phi)$, let $S(R)=\prod_{k=1}^n\phi_{i_k,d}'$ be the `shortest side' of it.

\begin{remark}\emph{ Let $R$ and $R'$ be two $n$-th cylinders of $\Lambda_\Phi$, and $\nu$ be
a Bernoulli measure of $\Lambda_\Phi$. Then $\nu(R \cap R')=0$ is always true.
See for instance \cite{Rao2019}.
}
\end{remark}

\begin{proof}[\textbf{Proof of Theorem \ref{thm:sponge}.}]
For $1\le j\le d$, define
\begin{equation}\label{alphaj}
\alpha_j:=\sum_{k=1}^j\beta_k.
\end{equation}
We shall prove by induction on $j$ that
for any cylinder $R$ of $\Lambda_j$ and any $\delta\le S(R)$,
\begin{equation}\label{eq:box-j}
{\mathcal{N}_\delta }(R) \asymp {\mu_j} (R) \delta^{-{\alpha_j}}.
\end{equation}

If $j=1$, $\Lambda_1$ is a one-dimensional self-similar set satisfying the
OSC, hence $\alpha_1=\beta_1=\dim_H \Lambda_1=\dim_B \Lambda_1$.
By Theorem \ref{thm:OSC},
$\mu_1$ is  a box-counting measure of $\Lambda_1$ and \eqref{eq:box-j} holds.

Now suppose  \eqref{eq:box-j} holds for $j=d-1$.
Let $R=\phi_{i_1\dots i_n}(\Lambda_\Phi)$ be a cylinder of $\Lambda_d=\Lambda_\Phi$.
Let us denote
$$
r_*:=\min\{\phi'_{\ba,d};~\ba\in\SD\}.
$$
We start with  the special case that $S(R)r_*< \delta\leq S(R).$
 A  crucial observation is that
$$\mathcal{N}_\delta(R)\asymp \mathcal{N}_\delta(\pi_{d-1}(R)).$$
Denote $\widetilde{R}=\pi_{d-1}(R)$, then $\widetilde{R}$ is a cylinder of $\Lambda_{d-1}$. By induction hypothesis, we have
$$
\mathcal{N}_\delta(\widetilde{R})\asymp \mu_{d-1}(\widetilde{R})\delta^{-\alpha_{d-1}}.
$$

For simplicity, let us denote $\phi_{i_k}=(f_{k,1},\dots, f_{k,d})$.
Then
\begin{equation}\label{eq:muC}
\mu_d(R)=\prod_{k=1}^n\prod_{j=1}^d(f'_{k,j})^{\beta_k}
\end{equation}
 and
 \begin{equation}\label{eq:SC}
 \delta\le S(R)=\prod_{k=1}^n f'_{k,d}< \delta/r_*.
 \end{equation}
It follows that
\begin{equation}\label{eq:C}
\begin{array}{rl}
\mathcal{N}_\delta(R)& \asymp \mu_{d-1}(\widetilde{R})\delta^{-\alpha_{d-1}}\\
            &=\left(\prod_{k=1}^n\prod_{j=1}^{d-1}(f'_{k,j})^{\beta_k}\right) \delta^{-\alpha_{d-1}}\\
           &=\left(\prod_{k=1}^n\prod_{j=1}^d(f'_{k,j})^{\beta_k}\right)  \left(\prod_{k=1}^n(f'_{k,d})^{-\beta_d}\right) \delta^{-\alpha_{d-1}}\\
           &\asymp \left(\prod_{k=1}^n\prod_{j=1}^d(f'_{k,j})^{\beta_k}\right) \delta^{-\beta_d} \delta^{-\alpha_{d-1}}  \quad(\text{by \eqref{eq:SC}})\\
           &=\mu_d(R) \delta^{-\alpha_{d}}, \quad(\text{by \eqref{alphaj} and \eqref{eq:muC}})
\end{array}
\end{equation}
which verifies \eqref{eq:box-j}.

Now we consider the general case that $\delta\le S(R)$.
Let  ${\mathcal V}_\delta$ be the collection of  cylinders $H$  such that
$S(H)\leq \delta < S(\hat H)$, where $\hat H$ is the ancestor of $H$; then
 $$\delta r_*<  S(H)\leq \delta.$$
 First,  by \eqref{eq:C} we have
 \begin{equation}\label{eq:H} \mathcal{N}_\delta(H)\asymp \mu_d(H) \delta^{-\alpha_{d}}\end{equation}
Clearly, ${\mathcal V}_\delta$ is a covering of $\Lambda_d$,
and the elements in ${\mathcal V}_\delta$ have disjoint interiors. Moreover, $R$ is a union of some elements of ${\mathcal V}_\delta$. These properties guarantee that
$$
\mathcal{N}_\delta(R)\asymp \sum_{H\in {\mathcal V}_\delta \text{ and } H\subset R} \mathcal{N}_\delta(H),
$$
which together with \eqref{eq:H} imply that
$$
\mathcal{N}_\delta(R)\asymp \sum_{H\in {\mathcal V}_\delta \text{ and } H\subset R} \mu_d(H)\delta^{-\alpha_{d}}
=\mu_d(R)\delta^{-\alpha_{d}}.
$$
So \eqref{eq:box-j} holds.
Consequently,  $\dim_B \Lambda_d=\alpha_d$.
 The theorem  is proved.
\end{proof}

\section{\textbf{Box-counting measures of Bara\'nski carpets}}\label{sec:Baranski}

 Let $r,s\geq 2$ be two integers.
Let $\{\Phi_{i,1}\}_{i=1}^r$ and $\{\Phi_{j,2}\}_{j=1}^s$ be the base IFS's such that
$\{\Phi_{k,1}[0,1); ~1\leq k\leq r\}$ and  $\{\Phi_{k,2}[0,1); ~1\leq k\leq s\}$ are
two disjoint families.
Let $\SD\subset \{1,\dots, r\}\times\{1,\dots,s\}$. For $\bd=(i,j)$, let $\Phi_{\bd}=(\Phi_{i,1},\Phi_{j,2})$. The attractor $\Lambda$ of the diagonal self-affine IFS $$\{\Phi_{\bd}\}_{\bd\in\SD}$$ is called a \emph{Bara\'nski carpet} (Bara\'nski \cite{Baranski07}).

Write $a_i=\Phi_{i,1}'$ and $b_j=\Phi_{j,2}'$. For $\bd=(i,j)$,  we denote $a_\bd=a_i$ and $b_\bd=b_j$; in other words, $\Phi_\bd'=(a_\bd,b_\bd)$.
Let $\pi_1(x,y)=x$ and $\pi_2(x,y)=y$ be the canonical projections.
For $I=i_1\dots i_k\in\SD^k,$ we denote
$$
a_{I}=a_{i_1}\cdot \cdots \cdot a_{i_k},
\quad b_{I}=b_{i_1}\cdot \cdots \cdot b_{i_k}, \text{ and }
s_{I}=\min\{a_{I},
b_{I}\}.
$$

Let $\SD^*=\bigcup_{k\geq 0} \SD^k$ with $\SD^0=\emptyset$. Set
$$
{\mathcal V}^A=\{I\in \SD^*;~ a_I\geq b_I\}, \quad  {\mathcal V}^B=\{I\in \SD^*;~ a_I< b_I\},
$$
and for    $0<\delta<1$,   set
$$
{\mathcal {V}_\delta}=\{I=i_1\dots i_k \in\SD^*;~
s_{I}\leq \delta< s_{i_1\dots i_{k-1}}\}
$$
be the collection of cylinders whose shorter sides are of `approximately equal size'.
Moreover, we define
$$
\mathcal {V}_\delta^A=\{{I\in{\mathcal {V}_\delta}:a_I\ge b_I}\},\quad
\mathcal {V}_\delta^B ={\mathcal {V}_\delta }\backslash{\mathcal {V}_\delta}^A.
$$

Let $\alpha_1$ and $\alpha_2$  be the positive real numbers satisfying
$$
\sum_{i\in \pi_1(\SD)} a_i^{\alpha_1}=1, \quad \sum_{\bd\in \SD} a_{\bd}^{\alpha_1} b_{\bd}^{\alpha_2}=1.
$$
Similarly, let $\beta_1$ and $\beta_2$  be the positive real numbers satisfying
$$
\sum_{j\in \pi_2(\SD)} b_j^{\beta_1}=1, \quad \sum_{\bd\in \SD}  b_{\bd}^{\beta_1}a_{\bd}^{\beta_2}=1.
$$
Set $\alpha=\alpha_1+\alpha_2$ and $\beta=\beta_1+\beta_2$. (Bara\'nski \cite{Baranski07} proved that
 $\dim_B \Lambda=\max\{\alpha,\beta\}.$)

Let $\nu_A$ be the Bernoulli measure on $\pi_1(\Lambda)$ with probability weight $(a_i^{\alpha_1})_{i\in \pi_1(\SD)}$,  and $\nu_B$  on $\pi_2(\Lambda)$ with probability weight $(b_j^{\beta_1})_{j\in \pi_2(\SD)}$.
Then $\nu_A$ and $\nu_B$ are box-counting measures since
they are Hausdorff measures of self-similar sets (Theorem \ref{thm:OSC}).

Let $\mu_A$ and $\mu_B$ be the Bernoulli measures on $\Lambda$ with probability weights
$$(a_\bd^{\alpha_1}b_\bd^{\alpha_2})_{\bd \in \SD} \text{ \ and \ }
(b_\bd^{\beta_1}a_\bd^{\beta_2})_{\bd \in \SD}$$
respectively.

\begin{lemma}\label{Bcover}
Let $\Lambda$ be a Bara\'nski carpet. Then
\begin{equation}\label{B_number}
\mathcal{N}_{b_I}(\Lambda_I)\asymp
 \mu_A(\Lambda_I)b_I^{-\alpha}   \text{ for } I\in {\mathcal V}^A,
\end{equation}
and
\begin{equation}\label{A_number}
\mathcal{N}_{a_I}(\Lambda_I)\asymp
 \mu_B(\Lambda_I)a_I^{-\beta}   \text{ for } I\in {\mathcal V}^B.
\end{equation}
\end{lemma}

\begin{proof} If $I\in {\mathcal V}^A$, then $a_I\geq b_I$. The crucial observation is that  $$\mathcal{N}_{b_I}(\Lambda_I)\asymp \mathcal{N}_{b_I}(\pi_1(\Lambda_I)).$$
Notice that $\pi_1(\Lambda)$ is a self-similar set satisfying the OSC,
 $\nu_A$  is  the canonical Bernoulli measure on $\pi_1(\Lambda)$, and $\dim_B\pi_1(\Lambda)=\alpha_1$.
By Theorem \ref{thm:OSC},  for $\delta\leq a_I$ we have
$$\mathcal{N}_\delta(\pi_1(\Lambda_I))\asymp   \frac{\nu_A(\pi_1(\Lambda_I))}{ \delta^{\alpha_1}}
=
\frac{a_I^{\alpha_1}}{ \delta^{\alpha_1}}.$$
Set $\delta=b_I$ in the above formula, we obtain \eqref{B_number}. Formula \eqref{A_number} can be obtained in the same manner.
\end{proof}


\begin{lemma}[Bara\'nski \cite{Baranski07}]\label{lem:BaBa}
Let $\Lambda$ be a Bara\'nski carpet. Then
\begin{equation*}
 \max\{\alpha, \beta\}\le{\alpha_1}+{\beta_1}.
\end{equation*}
\end{lemma}

\begin{proof}
Note that
$$
\sum_{\bd\in \SD} {a_{\bd}^{\alpha_1}b_\bd^{\beta_1}}
\le
\left(\sum\limits_{i\in \pi_1(\SD)}{a_i^{\alpha_1}}\right)\left(\sum\limits_{j\in \pi_2(\SD)}{b_j^{\beta_1}}\right)=1.
$$
So we have $\beta_1\geq \alpha_2=  \alpha-\alpha_1$ and $\alpha_1\geq \beta_2=  \beta-\beta_1$,
which imply the lemma.
\end{proof}

\begin{lemma}\label{lem:compare} Let $\Lambda$ be a Bara\'nski carpet. Then
$$
\mu_A(\Lambda_I) b_I^{-\alpha}\geq \mu_{B}(\Lambda_I)b_I^{-\beta}, \quad \text{ if }
I\in {\mathcal V}^A;
$$
$$
\mu_A(\Lambda_I) a_I^{-\alpha}\leq \mu_{B}(\Lambda_I)a_I^{-\beta}, \quad \text{ if }
I\in {\mathcal V}^B.
$$
\end{lemma}

\begin{proof}  Let $I\in \mathcal {V}^A$, then $a_I\geq b_I$. By Lemma \ref{lem:BaBa}, we have  $\alpha_1\geq \beta_2$,
so
$$
\frac{\mu_A(\Lambda_I) b_I^{-\alpha}}{\mu_{B}(\Lambda_I)b_I^{-\beta}}=
\left (\frac{a_I}{b_I}\right )^{\alpha_1-\beta_2}\geq 1,$$
 which implies the first inequality of the lemma.
 The second inequality can be obtained in the same manner.
\end{proof}

For a cylinder $R=\Lambda_I$, we  denote $S(R)=s_I$ to  be the length of its shorter side.
\medskip

\begin{proof}[\textbf{Proof of Theorem \ref{thm:Baranski}.}]
Without loss of generality, we assume that  $\alpha\geq \beta$. We will show that for any cylinder $R$,
$$\mathcal {N}_\delta(R)\asymp (\mu_A(R)+\mu_B(R))\delta^{-\alpha}, \quad \text{ for } \delta\leq S(R),$$
 in case of $\alpha=\beta$, and
$$\mathcal {N}_\delta(R)\asymp\mu_A(R)\delta^{-\alpha}$$
holds for $\delta\leq\delta_0(R)$ where
\begin{equation}\label{eq:delta_0}
\delta_0(R)=
\min\{S(R), (\mu_A(R)/\mu_B(R))^{\frac{1}{\alpha-\beta}}\}
\end{equation}
in case of $\alpha>\beta$.

Let $\delta\le S(R)$. Note that $R$ can be written as a finite disjoint union of cylinders in ${\mathcal V}_\delta$, that is, $R=\bigcup_{I\in\mathcal {V}_\delta}\Lambda_I$.
 Notice that for each  $I\in {\mathcal V}_\delta$, it holds that
$$r_*\delta<\min\{a_I, b_I\}\leq \delta,$$
where
$r_*:=\min(\{a_\bd; \bd\in \SD\}\cup\{b_\bd;~\bd\in \SD\}).$
By Lemma \ref{Bcover},  we have
\begin{equation}\label{total_num_1}
\begin{array}{rl}
\mathcal {N}_\delta(R) &\asymp  \sum_{I\in {\mathcal V}_\delta \text { and } \Lambda_I\subset R }\mathcal{N}_\delta(\Lambda_I) \\
&\asymp
 \sum_{I\in\mathcal {V}_\delta^A\text { and } \Lambda_I\subset R}\mu_A(\Lambda_I) \delta^{-\alpha}+ \sum_{I\in\mathcal {V}_\delta^B\text { and } \Lambda_I\subset R}\mu_B(\Lambda_I) \delta^{-\beta}.
\end{array}
\end{equation}
So there exists $M_1>0$ such that
\begin{equation}\label{total_num_2}
\mathcal {N}_\delta(R)\le M_1(\mu_A(R)\delta^{-\alpha}+\mu_B(R)\delta^{-\beta}).
\end{equation}

According to Lemma \ref{lem:compare}, we have
 \begin{equation*}
   \sum_{I\in\mathcal {V}_\delta^B\text { and } \Lambda_I\subset R}\mu_B(\Lambda_I) \delta^{-\beta}
 \geq    \sum_{I\in\mathcal {V}_\delta^B\text { and } \Lambda_I\subset R}\mu_A(\Lambda_I) \delta^{-\alpha},
 \end{equation*}
 so by \eqref{total_num_1} we obtain that there is a universal constant $M_2>0$ such that
$$
\mathcal {N}_\delta(R)\geq M_2^{-1}\mu_A(R)\delta^{-\alpha}.
$$
 By symmetry, we can show that
$\mathcal {N}_\delta(R)\geq M_3^{-1}\mu_B(R)\delta^{-\beta}$ for a universal constant $M_3>0$. It follows that
$$
\mathcal {N}_\delta(R)\geq\frac{\mu_A(R)\delta^{-\alpha}+\mu_B(R)\delta^{-\beta}}{M_2+M_3}.
$$
This together with \eqref{total_num_2} implies that
\begin{equation*}
\mathcal {N}_\delta(R)\asymp  \mu_A(R)\delta^{-\alpha}+\mu_B(R)\delta^{-\beta}.
\end{equation*}
In case of $\alpha=\beta$, this already gives us the desired relations.

Now suppose that $\alpha>\beta$. Choose $\delta_0$ as in \eqref{eq:delta_0}, then
 for $\delta\leq \delta_0$,
we have $\mu_A(R)\delta^{-\alpha}\geq \mu_B(R)\delta^{-\beta}$, so
 $$
\mathcal {N}_\delta(R)\asymp  \mu_A(R)\delta^{-\alpha}.
$$
 This completes the proof. As a by-product, we obtain that $\dim_B \Lambda=\alpha$.
\end{proof}

\section{\textbf{Box-counting measures of symbolic spaces}}\label{sec:symbolic}

 In this section, we consider the box-counting measures of several symbolic spaces related to
 self-affine sponges.

Let  $0<\xi<1$. Let $D_\xi$ be the metric on $\Z^\infty$ defined by
\begin{equation}\label{def:metric}
D_\xi(I, J)=\xi^{|I\wedge J|}, \quad I, J\in \Z^\infty,
\end{equation}
where $I\wedge J$ is the maximal common prefix of $I$ and $J$, and $|W|$ denotes the length of a finite word $W$.

\subsection{The first metric ($\lambda$-metric)}
\ \\
\indent Let $1>\xi_1>\xi_2>\cdots>\xi_d>0$   be  a sequence of real numbers.
The metric $\lambda$ on $(\Z^d)^\infty$ is defined  as the product metric
\begin{equation*}
((\Z^d)^\infty, \lambda)=(\Z^\infty,D_{\xi_1})\times \cdots \times (\Z^\infty, D_{\xi_d})
\end{equation*}
with $\lambda=\max\{D_{\xi_1},\dots,D_{\xi_d}\}$.

Let $\SD\in\mathbb{Z}^d$ be a finite set. Then $(\SD^\infty, \lambda)$
 is a  totally disconnected  metric space.
  This symbolic space
 has been used in many works related to Bedford-McMullen carpets, for example,
\cite{King95,JR11} on multi-fractal analysis,  and \cite{Paper-3, Yang22} on Lipschitz
classification.

 For $I\in \SD^k$, recall that $[I]$ is a $k$-th cylinder of
 $\SD^\infty$.
 Denote $N:=\#\SD$.  Let $\mu$ be the measure on $\SD^\infty$ such that for any cylinder $[I]$ of rank $k$,
  $\mu([I])=1/N^k,$
  and we call it the \emph{uniform Bernoulli measure} on $\SD^\infty$.

\begin{thm} The uniform Bernoulli measure is a box-counting measure of the   space  $(\SD^\infty, \lambda)$.
\end{thm}

\begin{proof} Similar to Section 4, we define $\beta_1,\dots, \beta_d$ inductively as following.
Let $\beta_1>0$ be the unique real number satisfying
$
  \#\left (\pi_1(\SD)\right )  (\xi_1)^{\beta_1}=1.
$
If $\beta_1,\dots, \beta_{j-1}$ are defined, we set $\beta_j>0$ to be the unique real number satisfying
$$
  \#\left (\pi_j(\SD)\right ) \prod_{k=1}^j (\xi_k)^{\beta_k}=1.
$$
Then by the same argument as the proof of Theorem \ref{thm:sponge}, one can show that $\beta=\sum_{j=1}^d \beta_j$
is the box dimension of $(\SD^\infty, \lambda)$, and the uniform Bernoulli measure
is a box-counting measure.
\end{proof}

\subsection{The second metric ($\rho$-metric)}
\ \\
\indent  To study the Lipschitz classification
of Bedford-McMullen carpets, \cite{Paper-3} introduced
  a quasi-metric $\rho$ on $\SD^\infty$ as following.
  Let $2\leq m<n$ be two integers.
  Assume that
$$
\SD \subset  \{0, 1,\dots, m-1\}\times \Z.
$$
For $I=(\bi,\bj), I'=(\bi',\bj')  \in \SD^\infty$, set
\begin{equation*}
 \rho(I, I')=\max\{r_{1/m}(\bi,\bi'), D_{1/n}(\bj,\bj')\},
\end{equation*}
where $D_{1/n}$ is defined by \eqref{def:metric}
 and
$$r_{1/m}(\bi,\bi')=\left |\sum_{k=1}^\infty\frac{i_k}{m^k}-\sum_{k=1}^\infty\frac{i_k'}{m^k}\right |.$$
 If for any two distinct points $I=(\bi,\bj), I'=(\bi',\bj')\in \SD^\infty$,
it holds that  $\rho(I, I')\neq 0$, then we say $\SD$ satisfies the \emph{non-overlapping condition}; in this case $(\SD^\infty, \rho)$
is a metric space.

\begin{thm} Suppose $\SD$ satisfies the non-overlapping condition. Then the uniform Bernoulli measure is a box-counting measure of the   space $(\SD^\infty, \rho)$.
\end{thm}
\begin{proof}
Let $s=\#(\pi_1(\SD))$,  $\beta_1=\log s/\log m$,
 and
$\beta_2= \log(N/s)/\log n$  where $N=\#\SD$.
By the same argument as the proof of Theorem \ref{thm:sponge}, one can show that $\beta=\beta_1+\beta_2$
is the box dimension of $(\SD^\infty, \rho)$, and the uniform Bernoulli measure
is a box-counting measure.
\end{proof}

\section{\textbf{Relation to Minkowski content}}\label{sec:Minkowski}

 For $A\subset\mathbb{R}^d$, and    $0\le \beta \le d$, the \emph{$\beta$-dimensional upper Minkowski content} is
$$
\mathcal{M}^{*\beta}(A)=\limsup_{\delta\to 0}\frac{{\cal L}^d([A]_\delta)}{\delta^{d-\beta}}
$$
where ${\cal L}^d$ is the $d$-dimensional Lebesgue measure.
By taking lower limit instead of upper limit, we define
 the \emph{$\beta$-dimensional lower Minkowski content}
 $ \mathcal{M}_*^{\beta}(A)$.
(See for instance, $\S$3.1 of \cite{Fal90}.)

The following lemma  provides a necessary condition for the existence of  the box-counting measure.

\begin{lemma}\label{lem:content} Let $X$ be a compact subset of $\R^d$ with $\beta=\dim_B X$. If $\mu$ is a box-counting measure of $X$,
then there exists a constant $M>0$ such that
for any $\delta$-connected component $R$ of $X$, it holds that
\begin{equation*}
M^{-1} \mu(R)\leq  {\mathcal M}_*^{\beta}(R)\leq    {\mathcal M}^{*\beta}(R)\leq M \mu(R).
\end{equation*}
\end{lemma}

\begin{proof} For $A\subset X$, define
$$
g_\delta(A)= \frac{{\cal L}^d([A]_\delta)}{\delta^d}.
$$
It is easy to show that $g_\delta(A)\asymp \mathcal{N}_\delta(A).$
Hence $g_\delta(R)\asymp \mu(R)\delta^{-\beta},$
which proves the lemma.
\end{proof}

\end{document}